\documentclass[a4paper,twoside]{article}

\usepackage{geometry}
\usepackage{calc}
\usepackage{amssymb}
\usepackage{amstext}
\usepackage{amsmath}
\usepackage{amsthm}
\usepackage{graphicx}
\usepackage{booktabs}
\usepackage{array}
\usepackage{caption}
\usepackage{float}
\usepackage{multirow}
\usepackage{setspace}
\usepackage{pifont}
\usepackage{breqn}
\usepackage{color}
\usepackage{algorithm}
\usepackage{algorithmic}
\usepackage{tikz}
\usepackage{pgfgantt}
\usepackage{wrapfig}
\usepackage{adjustbox}
\usepackage[affil-it]{authblk}

\usepackage{natbib}

\newcommand{\ignore}[1]{}

\providecommand{\keywords}[1]{\textbf{\textit{Keywords:}} #1}
\providecommand{\abstractNew}[1]{\textbf{\textit{Abstract:}} #1}

\theoremstyle{definition}

\begin{document}

\title{\bf Bi-objective risk-averse facility location using a subset-based representation of the conditional value-at-risk}

\author[1]{Najmesadat Nazemi\thanks{najmesadat.nazemi@jku.at}}
\author[1]{Sophie N. Parragh\thanks{sophie.parragh@jku.at}}
\author[2]{Walter J. Gutjahr\thanks{walter.gutjahr@univie.ac.at}}

\affil[1]{Institute of Production and Logistics Management, Johannes Kepler University Linz, Austria}
\affil[2]{Department of Statistics and Operations Research,  University of Vienna, Austria}

\date{}

\maketitle

\keywords{Risk-averse, Conditional value-at-risk, Cutting-plane, Bi-objective optimization\\}

\abstractNew{For many real-world decision-making problems subject to uncertainty, it may be essential to deal with multiple and often conflicting objectives while taking the decision-makers' risk preferences into account. Conditional value-at-risk (CVaR) is a widely applied risk measure to address risk-averseness of the decision-makers.
In this paper, we use the subset-based polyhedral representation of the CVaR to reformulate the bi-objective two-stage stochastic facility location problem presented in \cite{nazemi2020bi}.
We propose \textcolor{black}{an approximate} cutting-plane method to deal with this more computationally challenging subset-based formulation\textcolor{black}{. Then, the cutting plane method is embedded into the $\epsilon$-constraint method, the balanced-box method, and a recently developed matheuristic method to address the bi-objective nature of the problem}.
Our computational results show the effectiveness of the proposed method. Finally, we discuss how incorporating an approximation of the subset-based polyhedral formulation affects the obtained solutions.\\[1ex]

{\bf  \textit{Published in}} Proceedings of the 11th International Conference on Operations Research and Enterprise Systems - ICORES, pages  77-85, 2022, doi:10.5220/0010914900003117
}

\section{\uppercase{Introduction}}
\label{sec:intro}

On the one hand, many real-world problems have multiple conflicting objectives, e.g., cost versus customer satisfaction or profitability versus environmental concerns. On the other hand, it is important to incorporate uncertainty about the considered parameters into the decision-making process.
Optimization under uncertainty and multi-objective optimization have evolved mostly separately over the past decades.
However, it is desirable to develop optimization techniques to deal with these two domains simultaneously \cite{gutjahr2016stochastic}.

Due to its sound theoretical background and its proven algorithmic performance, stochastic programming is a widely applied approach in computational optimization for handling data uncertainty 
(see, e.g., \cite{birge2011introduction} and \cite{shapiro2014lectures}, and references therein). The standard two-stage stochastic programming problem formulation considers the expected value as a performance measure.  
It assumes 
a risk-neutral attitude of the decision-maker and performs well in the context of repetitive decision-making problems \citep{birge2011introduction}. However, decision-makers may 
have other  
risk preferences, especially in non-repetitive applications, e.g., in large-scale financial investments or in the management of high-impact disasters.
In the case where more robust solutions are of interest, one of the recently most widely recommended risk measures is the Conditional Value at Risk (CVaR). It can be adjusted to the actual degree of risk-averseness of the decision-maker by specifying 
a confidence level $\alpha$ ($\alpha \in [0,1)$) (CVaR($\alpha$)). The closer the value of $\alpha$ to 0, the more risk-averse is the decision-maker.
CVaR has been applied in different streams of the optimization literature, such as, e.g., financial engineering (e.g., \citet{krokhmal2002portfolio}, \citet{roman2007mean}, \textcolor{black}{\citet{huang2021multi}, \citet{zheng2021novel}}) disaster management (e.g., \citet{noyan2012risk}, \citet{elcci2018chance}, \citet{nazemi2020bi}), water management (e.g., \citet{zhang2016decomposition}, \citet{simic2019interval}, \citet{chen2021coupling}), and sustainable supply chain (e.g., \citet{rahimi2018sustainable}, \citet{rahimi2019stochastic}). \textcolor{black}{We refer to \citet{filippi2020conditional} for a recent survey on application of CVaR to several optimization domains different from financial engineering.}
\textcolor{black}{\citet{huang2021multi} propose a two-stage distributionally robust model including CVaR as a risk measure for two different practical applications, a multi-product assembly problem and a portfolio selection problem. They also extend the model to a multi-stage one.  Two decomposition algorithms based on the cutting-plane method are developed to solve the models. Finally, they conclude that the proposed models show a more robust performance compared to the risk-neutral model.}
\cite{noyan2012risk} incorporates the classical CVaR representation into the objective function of a two-stage mean-risk stochastic problem in disaster management where the demand and the damage level of the transportation network are uncertain parameters. She develops decomposition algorithms to deal with the computationally challenging nature of the problem. In line with this study, \cite{elcci2018chance} consider a chance-constrained two-stage mean-risk stochastic relief network design model. They also develop a Benders decomposition-based algorithm to solve formulations with alternative representations of the CVaR. \cite{zhang2016decomposition} present a risk-averse multi-stage stochastic water allocation problem. A mean-CVaR objective is considered in each stage of the problem. Finally, they propose a nested L-shaped method to solve the problem.
In addition to the aforementioned works that address single-objective problems, there are also some studies in the literature that consider multi-objective optimization under uncertainty incorporating CVaR as the risk measure into two-stage stochastic programming.
\textcolor{black}{\cite{zheng2021novel} consider a bi-objective portfolio optimization model with a focus on minimizing CVaR as a risk measure and maximizing the expected return rate. They also incorporate transaction costs into the objective functions. A NSGA-II metaheuristic method \citep{deb2002fast} based on sparsity strategy \citep{zitzler2001spea2} (SMP-NSGAII) is developed to deal with the bi-objective model.}
\cite{rahimi2019stochastic} propose a risk-averse sustainable multi-objective mixed-integer non-linear model to design a supply chain network under uncertainty by incorporating CVaR into their two-stage stochastic model. The objectives aim at minimizing the design cost and maximizing the profit. 
\cite{nazemi2020bi} incorporate a risk-neutral, a risk-averse, and the CVaR measure into a bi-objective two-stage facility location model in a disaster management context to analyze a wide range of risk preferences, including risk-neutral and worst-case approaches. They integrate different uncertain two-stage models into two well-known exact multi-objective frameworks, namely the $\epsilon$-constraint and the balanced-box methods. Additionally, they also develop a matheuristic approach and 
they analyze and evaluate the combination of different uncertainty representations and multi-objective frameworks.
In most of the aforementioned studies, the classical representation of CVaR has been employed to model risk-averseness.
However, as mentioned earlier, alternative variants of CVaR are also proposed in the literature. In this paper, we focus on the existing two-stage bi-objective model presented in \cite{nazemi2020bi} and extend this work 
by reformulating their two-stage risk-averse model using  
a subset-based polyhedral representation of the CVaR. 
To address this alternative formulation, we propose a scenario cutting-plane algorithm. 
We conduct a computational analysis on the test instances
used by \cite{nazemi2020bi} to illustrate how the subset-based variant performs in comparison to the classical formulation on their model.

The remainder of this paper is organized as follows. In Section~\ref{sec:CVaR}, we describe the existing bi-objective facility location problem under uncertainty presented by \cite{nazemi2020bi} along with its corresponding subset-based polyhedral mathematical reformulation. In Section \ref{sec:method}, we summarize the developed solution approach to solve the problem. In Section \ref{result} we discuss the computational results. Finally, we present the conclusion and address potential future work in Section \ref{sec:conclusion}.

\section{\uppercase{Two-stage Risk-averse Stochastic Optimization}}
\label{sec:CVaR}

The problem addressed in \cite{nazemi2020bi} is a bi-objective two-stage location-allocation model under demand uncertainty motivated by last-mile network design in a slow-onset disaster context. 
A finite discrete set of demand scenarios ($s\in S$, $S=\{1,..., N\}$) with equal probabilities ($\frac{1}{N}$) is used to incorporate stochastic information.
The problem aims at finding the best location to position temporary local distribution centers (LDC) in the response phase of a disaster, such that the affected people can walk to these centers to receive their relief aids. The model tries to find the trade-off between two conflicting objectives, i.e., minimization of operational cost on opening LDCs and maximization of the coverage or in other words, minimization of the amount of uncovered demand. We use the same notation as introduced in \cite{nazemi2020bi} to describe the model. The problem is formulated on a network, where the sets of nodes representing the demand and potential LDC locations are denoted by $I$ and $J$, respectively, assuming that  
$J\subseteq I$. The assumption is that the affected people in each demand node ($i\in I$) will walk only to an LDC 
if their distance ($d_{ij}$) from an opened LDC is less than a certain distance threshold ($d_{max}$). Uncertainty is assumed on the amount of demand at each node ($q_i^{(s)}$).
The first-stage decision (here-and-now) 
concerns the selection of sites to open LDCs. Binary variables $y_j \in \{0,1\}$ indicate whether or not an LDC is opened at node $j \in J$. Each LDC
is assumed to have a limited capacity  $c_j$ and a given fixed cost $\gamma_j$. Decision variables $x_{ij}^{(s)} \in \{0,1\}$ and $u_j^{(s)}$ represent the second-stage decisions (wait-and-see) that determine the allocation of demand nodes to the opened LDCs and the amount of relief items delivered to each of them, respectively. The second-stage decisions are determined based on first-stage values and realized uncertain demand information. We next elaborate on our choice of the alternative subset-based polyhedral representation of CVaR incorporated into the second stage of the stochastic model.

\subsection{Two-stage Subset-based CVaR (SSCVaR)}

\cite{fabian2008handling} obtains a subset-based polyhedral representation of the classical CVaR for the special case of scenarios with equal probabilities where $\alpha=1-\frac{k}{N}$; $k$ is the cardinality of any subset of scenarios. The cardinality of all scenarios is equal to $N$ (i.e., $\bar{S}\subseteq S , |\bar{S}|=k$ and $|S| = N$):

\begin{equation}\label{eq:CVaR}
\text{CVaR}_{1-\frac{k}{N}}(X)=\frac{1}{k} \max_{\bar{S}\subseteq S : |\bar{S}|=k} \sum_{s \in \bar{S}}X^{(s)}
\end{equation}

where $X^{(s)}$ is the value of the random variable $X$ in scenario $s$. Equation (1) leads to the following subset-based polyhedral representation of CVaR \citep{fabian2008handling}:

\begin{equation}\label{eq:subsetCVaR}
\text{CVaR}_{1-\frac{k}{N}}(X)=\{\min \varrho: \varrho \geq \frac{1}{k} \sum_{s \in \bar{S}}X^{(s)},  \forall \bar{S}\subseteq S: |\bar{S}|=k\}
\end{equation}

We utilize the subset-based representation \eqref{eq:subsetCVaR} to obtain the MIP formulation of the two-stage model (\textbf{MB}):

\begin{equation}\label{eq3}
    \text{min}\hspace{2mm}\sum_{j \in V} \gamma_{j}y_{j}
\end{equation}

\begin{equation}\label{eq4}
    \text{min}\hspace{2mm}\text{CVaR}_\alpha(\sum_{i \in I} q_i^{(s)}-\sum_{j \in J} u_j^{(s)}) = \varrho
\end{equation}

	\begin{equation}\label{eq5}
	\text{s.t.}\hspace{1mm}\varrho \geq \frac{1}{k} \sum_{s \in \bar{S}}(\sum_{i \in I} q_i^{(s)}-\sum_{j \in J} u_j^{(s)})\hspace{1.5mm}\forall \bar{S}\subseteq S: |\bar{S}|=k
	\end{equation}
	
	\begin{equation}\label{eq6}
	   \hspace{-1mm} \sum_{j \in J} \psi(d_{ij})x_{ij}^{(s)}\leq 1 \hspace{8mm}\forall i \in I, s \in S
	\end{equation}
	
	\begin{equation}\label{eq7}
	   u_j^{(s)}\leq c_j y_j\hspace{17mm}\forall j \in J, s \in S  
	\end{equation}
	
	\begin{equation}\label{eq8}
	   u_j^{(s)}\leq \sum_{i \in I} q_i^{(s)}\psi(d_{ij})x_{ij}^{(s)}\hspace{2mm}\forall j \in J, s \in S 
	\end{equation}
	
	\begin{equation}\label{eq9}
	   \hspace{1mm}x_{ij}^{(s)}\leq y_j \hspace{13mm}\forall i\in I, \forall j \in J, s \in S 
	\end{equation}
	
	\begin{equation}\label{eq10}
	  \hspace{-16mm} y_j\in \{0,1\}\hspace{10mm}\forall j \in J 
	\end{equation}
	
	\begin{equation}\label{eq11}
	 \hspace{6mm} x_{ij}^{(s)}\in \{0,1\}\hspace{8mm}\forall i\in I, \forall j \in J, s \in S
	\end{equation}
	
	\begin{equation}\label{eq12}
	  \hspace{-9mm} u_j^{(s)}\in \mathbb{Z^+}\hspace{13mm}\forall j \in J, s \in S 
	\end{equation}
	
	\begin{equation}\label{eq13}
	  \hspace{-10mm}\varrho \geq 0  
	\end{equation}

The first objective \eqref{eq3} minimizes total opening costs \textcolor{black}{of LDCs}. The second objective \eqref{eq4} minimizes the CVaR of the uncovered demand of the affected people. \textcolor{black}{The auxiliary variable $\varrho$ represents the solution value of the r.h.s. of \eqref{eq:subsetCVaR}, which coincides with the CVaR as given by \eqref{eq:CVaR}}. Constraints \eqref{eq5} make sure that $\varrho$ is set to the correct value.
Constraints \eqref{eq6} ensure that the demand at node $i$ can only be covered at most once by an opened LDC if it is located within the coverage threshold of LDC $j$ ($\psi({d_{ij}})$=1, if $d_{ij} \leq d_{max}$ and 0 otherwise). Constraints \eqref{eq7} are the capacity constraints, one for each LDC $j$. Constraints \eqref{eq8} guarantee that the amount of demand considered to be covered by LDC $j$ is at most the actual demand assigned to this node \textcolor{black}{within its coverage threshold}. Constraints \eqref{eq9} guarantee that a demand node can only be assigned to an opened LDC. Constraints \eqref{eq10}-\eqref{eq13} give the domains of the decision variables.

\section{\uppercase{Solution approach}}
\label{sec:method}

In this section, we describe the designed cutting-plane (approximation) method to deal with the potentially large number of scenario subsets in the model. 
To address the bi-objective nature of the problem, we embed the  cutting-plane method into the multi-objective frameworks applied in \cite{nazemi2020bi}, namely the $\epsilon$-constraint (e) \citep{laumanns2006efficient}, the balanced-box \citep{boland2015criterion} (BB) and \textcolor{black}{a specific matheuristic (Mat) developed in \cite{nazemi2020bi}}. \textcolor{black}{The matheuristic method embeds the bi-objective generalization of two heuristic approaches, namely local branching \citep{fischetti2003local} and relaxation induced neighborhood search (RINS) \citep{danna2005exploring}, into the $\epsilon$-constraint method. The RINS is a variable fixing scheme employed when moving along the Pareto frontier from one solution to another in the $\epsilon$-constraint framework. Additionally, a local branching approach is applied for the solutions, for which only fixing the variables does not solve the problem to optimality within a certain time limit.}
In a bi-objective setting, it is usually not possible to find a solution which simultaneously optimizes all of the considered objectives, but a set of optimal or efficient trade-off solutions exists, which are incomparable among each other and together dominate all other feasible solutions. The image of an efficient solution in objective space is called non-dominated point (NDP) and all NDPs together form the \emph{Pareto frontier}. 
For further background on multi-objective optimization, we refer to \cite{ehrgott2005multicriteria}.

\subsection{Scenario Subset-based Approximation}
Formulation (MB) contains constraints corresponding to subsets of $S$. The size of this set of constraints could grow fast, depending on the cardinality of the subsets ($k$) and the set of scenarios $S$. 

As mentioned in Section~\ref{sec:intro}, different decomposition-based techniques have been proposed to solve two-stage stochastic programming models with a single objective (e.g., \cite{shapiro2014lectures}). \cite{elcci2018chance} propose a Benders decomposition-based framework for the subset-based variant of CVaR with a mean-risk objective function. In the single objective two-stage stochastic model, the model is decomposed by scenario once the first-stage variables are fixed.

In our scenario subset cutting-plane approach, we iteratively solve the so-called single objective problem which relaxes the subset-based constraints \eqref{eq5}, and dynamically generates them using a delayed cut generation algorithm. 
To be more precise, the single objective problem contains the objective function \eqref{eq4} as the main objective, where the value of the other objective \eqref{eq3} is fixed.

Given an optimal solution $( \hat{y} , \hat{\varrho} ,\hat{u} )$ of the relaxed model in each iteration, we check whether there is any violated constraint of the type \eqref{eq5}, and identify a subset $\hat{S}$ if there is a violation. To solve this so-called \textit{separation} problem,
a few simple calculations (steps 3-8 of Algorithm \ref{alg4}) are required. At each iteration, we have the exact optimal second stage objective values associated with the current $( \hat{y} , \hat{\varrho} ,\hat{u} )$ vector; in particular, these objective values can be expressed as $\varrho^s=\frac{1}{k} (\sum_{i \in I} q_i^{(s)}-\sum_{j \in J} \hat{u}_j^{(s)})$ for each $s \in S $. Let $S^*= \{ s \in S :\sum_{s \in S^*} \varrho^s > \hat{\varrho} \} $. If $S^*= \{ s \in S : \sum_{s \in S^*} \varrho^s  \leq \hat{\varrho} \} $ then the candidate solution is labeled as the new incumbent. Otherwise, the identified violated constraint of the form \eqref{eq5} associated with $S^*$ is introduced as a cutting-plane. The pseudo-code of the delayed cut generation algorithm is presented in Algorithm \ref{alg4}.

\paragraph{Initial cut} We propose a heuristic algorithm in order to initialize the relaxed model. This procedure is based on the following steps:
\begin{enumerate}
	
	\item Rank the scenarios in each demand point based on their corresponding demand value ($q_{j}^s$) in descending order. 
	\item Sum the ranks for each scenario over $I$, the set of demand points. 
	\item Sort the scenarios based on their cumulative total rank in descending order.
	\item Select the first $k$ scenarios of the sorted set.
\end{enumerate}

\floatname{algorithm}{Algorithm}
\begin{algorithm}
	\caption{Delayed cut generation}
	\begin{algorithmic}[1]
		\STATE {initialize the relaxed model with the initial cut, $t \gets$ 1}\\
		\WHILE {$ \sum_{s \in \hat{S}} \varrho^s > \hat{\varrho}$}
		\FOR {each $s\in S$}
		\STATE {$\varrho^s \gets \frac{1}{k}{(\sum_{i \in V} q_i^{(s)}-\sum_{j \in V} \hat{u}_j^{(s)})}$}
		\STATE (sort $\varrho^s$ in descending order)
		\ENDFOR
		\STATE {$S^* \gets$ the first $k$ values (largest) of the sorted $\varrho^s$ values}
		\STATE $\hat{S} \gets S^*$ 
		\STATE $t \gets t+1$
		\ENDWHILE
	\end{algorithmic}
	\label{alg4}
\end{algorithm}

In our bi-objective setting, due to the presence of integer variables not only in the first stage, but also in the second stage, using cutting-plane technique for every single point along the Pareto frontier causes computationally no benefit. However, preliminary tests showed that, employing the cutting-plane algorithm in order to calculate the first extreme point of the Pareto frontier and, consequently, fixing the generated subset cuts to the model, leads to a high-quality approximation of the entire Pareto frontier. \textcolor{black}{We refer to this approach towards solving model MB by $\overline{\text{MB}}$ in the following.}

\section{\uppercase{Computational study}}
\label{result}
In this section, we present computational results for the instances used by \cite{nazemi2020bi} to show the performance of the proposed method with respect to the alternative subset-based polyhedral reformulation of the model in comparison to the classical representation of CVaR.

All experiments are executed on a single thread cluster where each node consists of two Intel Xeon X5570 CPUs at 2.93 GHz and 8 cores, with 48GB RAM. The model is implemented in C++ using CPLEX 12.9. As in \cite{nazemi2020bi} a time limit (TL) of 7200 s is considered. We note that we use the generic callback feature of CPLEX in order to generate the cuts in the MB model. Unlike lazy constraint callbacks, which prevent CPLEX from utilizing “dynamic search”, generic callbacks allow CPLEX to choose among the options of dynamic search and classical branch-and-cut as it suits. 

\subsection{Test Instances}
\label{test instances}
Here, we briefly explain the key points of the test instances and refer the reader to \cite{nazemi2020bi} for further details. They are derived from a real-world data set from a drought case presented in \cite{tricoire2012bi} for the region of Thies in western Senegal. The region is divided into 32 rural areas where each has between 9 and 31 villages. In total, the region contains 500 nodes. Opening costs for LDCs are assumed to be identical for all locations (5000 cost units). In addition, it is assumed that the affected people in a demand node walk to the closest LDC if the distance is less than 6km. For each instance, first, 10 sample scenarios ($|S|= 10$) are generated to cope with demand uncertainty. Then, in order to compare the algorithms for a larger number of scenarios, samples with sizes of $|S|$= 100, 500, and 1000 are generated.

\subsection{Solution Quality Indicators}
In order to assess the performance of heuristic methods for bi- and multi-objective optimization, several different quality indicators have been proposed in the literature (see, \cite{zitzler2003performance}).
In this study, instead of the hypervolume indicator ($I_H$) used by \cite{nazemi2020bi}, we employ the hypervolume gap (gH) and the multiplicative $\epsilon$-indicators.

For a bi-objective framework with minimization objective functions, assume $X$ is the set of feasible solutions in decision space and $Z=f(X)$ is the image of this set in objective space. We assume $A \subseteq Z$ is an approximation set of a Pareto frontier, and $R \subseteq Z$ is a reference set.

\textcolor{black}{The hypervolume gap (gH) definition is as follows:}

\begin{align}
  \text{gH\%}= \frac{I_H(R)-I_H(A)}{I_H(R)}\cdot100  
\end{align}

\textcolor{black}{The closer the value of gH\% to 0, the higher is the quality of the approximated Pareto frontier.}

The multiplicative epsilon indicator ($I_\epsilon$) \citep{zitzler2003performance} gives the minimum factor $\epsilon$ that the approximation set (A) in the objective space has to be multiplied with, such that it dominates the reference set (R). The closer the value of $I_\epsilon$ to 1, the better is the approximated Pareto frontier. The $\epsilon$ value calculated here can be interpreted as the gap between $A$ and $R$ in the objective space.

\subsection{Results}
We compare all combinations of multi-objective criterion space approaches used in \cite{nazemi2020bi} (e, BB, and Mat) with two linear counterparts of the classical (MA) and \textcolor{black}{the approximation of subset-based polyhedral ($\overline{\text{MB}}$)} representations of the two-stage CVaR model for different levels of risk. \{e-MA, BB-MA, Mat-MA, e-$\overline{\text{MB}}$, BB-$\overline{\text{MB}}$, Mat-$\overline{\text{MB}}$\} denotes the set of all methods evaluated in this study. 

The run time for the two different representations of the CVaR approach, MA, and $\overline{\text{MB}}$, is reported in Table~\ref{tab4}. MA is solved to optimality, whereas $\overline{\text{MB}}$ is the MB which is solved using the explained cutting-plane technique to approximate its Pareto frontier. \textcolor{black}{It is worth mentioning that the cutting-plane method solves the first extreme point of the Pareto frontier to optimality and all other non-dominated solutions are computed without generating additional cuts.} Results for different levels of  $\alpha$ are reported: the two special cases identical with the expected value (risk-neutral) and the worst-case (risk-averse) point of view, and one "middle" level and their corresponding $k$-values ($\alpha= 1- \frac{k}{N}$). The results in this table show that the combination of the BB with $\overline{\text{MB}}$ 
\textcolor{black}{solves the instances slightly faster} than its combination with MA. \textcolor{black}{In the following, we compare the quality of the obtained solutions using MA and $\overline{\text{MB}}$.}

% Table generated by Excel2LaTeX from sheet 'Sheet3'
\begin{table*}[!h]
	\centering
	\caption{Run time comparison of two representations of CVaR, MA, and $\overline{\text{MB}}$, for different $\alpha$-level/$k$-value for instances of size 21-500 with sample size 10}
		\small\addtolength{\tabcolsep}{-3pt}
	\makebox[\linewidth]{
		\scalebox{0.85}{
	\begin{tabular}{c|cccccccccccccc}
		\multicolumn{1}{r}{} & \multicolumn{9}{l}{$\alpha$/$k$ value}                            &       &       &       &       &  \\
		\cmidrule{2-15}    \multicolumn{1}{c}{} & \multicolumn{4}{c}{\textbf{0/10}} &       & \multicolumn{4}{c}{\textbf{0.7/3}} &       & \multicolumn{4}{c}{\textbf{0.9/1}} \\
		\cmidrule{1-5}\cmidrule{7-10}\cmidrule{12-15}    \multicolumn{1}{c}{\textit{\#Node}} & \textbf{e-MA} & \textbf{BB-MA} & \textbf{e-$\overline{\text{MB}}$} & \textbf{BB-$\overline{\text{MB}}$} &       & \textbf{e-MA} & \textbf{BB-MA} & \textbf{e-$\overline{\text{MB}}$} & \textbf{BB-$\overline{\text{MB}}$} &       & \textbf{e-MA} & \textbf{BB-MA} & \textbf{e-$\overline{\text{MB}}$} & \textbf{BB-$\overline{\text{MB}}$} \\
		\cmidrule{1-5}\cmidrule{7-10}\cmidrule{12-15}    \textbf{21} & 1     & 2     & 1     & 2     &       & 2     & 1     & 1     & 2     &       & 2     & 2     & 1     & 1 \\
		\textbf{44} & 16    & 22    & 15    & 19    &       & 17    & 23    & 13    & 14    &       & 20    & 26    & 10    & 17 \\
		\textbf{56} & 24    & 31    & 19    & 23    &       & 18    & 22    & 15    & 23    &       & 19    & 33    & 14    & 16 \\
		\textbf{72} & 36    & 51    & 29    & 35    &       & 45    & 48    & 30    & 41    &       & 64    & 103   & 36    & 40 \\
		\textbf{90} & 62    & 110   & 59    & 80    &       & 73    & 78    & 59    & 75    &       & 99    & 179   & 64    & 73 \\
		\textbf{106} & 123   & 229   & 118   & 133   &       & 151   & 161   & 111   & 132   &       & 210   & 380   & 141   & 175 \\
		\textbf{120} & 143   & 241   & 135   & 163   &       & 144   & 158   & 121   & 160   &       & 138   & 262   & 100   & 115 \\
		\textbf{163} & 445   & 818   & 392   & 552   &       & 417   & 506   & 369   & 503   &       & 636   & 1365  & TL    & 596 \\
		\textbf{182} & 577   & 1550  & 481   & 738   &       & 932   & 1053  & 701   & 776   &       & 1326  & 2548  & 1010  & 1192 \\
		\textbf{203} & 956   & 2078  & 879   & 1208  &       & 616   & 699   & 552   & 716   &       & 515   & TL    & TL    & 558 \\
		\textbf{254} & 6103  & TL    & TL    & TL    &       & TL    & TL    & TL    & TL    &       & TL    & TL    & TL    & TL \\
		\textbf{264} & 3903  & 6800  & TL    & 2449  &       & 6453  & TL    & TL    & TL    &       & TL    & TL    & TL    & 4007 \\
		\textbf{275} & 7133  & TL    & 6585  & TL    &       & TL    & TL    & TL    & TL    &       & TL    & TL    & TL    & TL \\
		\textbf{295} & TL    & TL    & TL    & TL    &       & TL    & TL    & TL    & TL    &       & TL    & TL    & TL    & TL \\
		\textbf{326} & TL    & TL    & TL    & TL    &       & TL    & TL    & TL    & TL    &       & TL    & TL    & TL    & TL \\
		\textbf{355} & TL    & TL    & TL    & TL    &       & TL    & TL    & TL    & TL    &       & TL    & TL    & TL    & TL \\
		\textbf{388} & TL    & TL    & TL    & TL    &       & TL    & TL    & TL    & TL    &       & TL    & TL    & TL    & TL \\
		\textbf{410} & TL    & TL    & TL    & TL    &       & TL    & TL    & TL    & TL    &       & TL    & TL    & TL    & TL \\
		\textbf{436} & TL    & TL    & TL    & TL    &       & TL    & TL    & TL    & TL    &       & TL    & TL    & TL    & TL \\
		\textbf{449} & TL    & TL    & TL    & TL    &       & TL    & TL    & TL    & TL    &       & TL    & TL    & TL    & TL \\
		\textbf{472} & TL    & TL    & TL    & TL    &       & TL    & TL    & TL    & TL    &       & TL    & TL    & TL    & TL \\
		\textbf{482} & TL    & TL    & TL    & TL    &       & TL    & TL    & TL    & TL    &       & TL    & TL    & TL    & TL \\
		\textbf{500} & TL    & TL    & TL    & TL    &       & TL    & TL    & TL    & TL    &       & TL    & TL    & TL    & TL \\
		\bottomrule
		\bottomrule
	\end{tabular}%
}
}
	\label{tab4}%
\end{table*}%

After analyzing the two exact multi-objective techniques (e and BB), we apply the Mat method in combination with the MA and $\overline{\text{MB}}$ models to solve the test instances. As shown in  
Table~\ref{tab5}, 
the $\overline{\text{MB}}$-Mat method  
is slightly faster than MA-Mat. Additionally, we solve a few instances with a larger number of scenarios. We report the run time and the number of found non-dominated points (NDP) in Table~\ref{tab6}. The results show that the approximately solved subset-based representation of the CVaR method ($\overline{\text{MB}}$) pays off when the number of scenarios increases.

% Table generated by Excel2LaTeX from sheet 'Sheet15'
\begin{table*}[!h]
	\centering
	\caption{Run time comparison of a combination of $\epsilon$-constraint method (e), and the proposed matheuristic (Mat) method with MA and $\overline{\text{MB}}$ for different $\alpha$/$k$-values for test instances of size 21-500 with sample size of 10 scenarios.}
		\small\addtolength{\tabcolsep}{-3pt}
	\makebox[\linewidth]{
		\scalebox{0.85}{
	\begin{tabular}{c|cccccccccccccc}
		\multicolumn{1}{r}{} & \multicolumn{8}{l}{$\alpha$/$k$ value}                             &       &       &       &       &       &  \\
		\cmidrule{2-15}    \multicolumn{1}{c}{} & \multicolumn{3}{c}{\textbf{0/10}} &       &       & \multicolumn{3}{c}{\textbf{0.7/3}} &       &       & \multicolumn{3}{c}{\textbf{0.9/1}} &  \\
		\cmidrule{1-5}\cmidrule{7-10}\cmidrule{12-15}    \multicolumn{1}{c}{\textit{\#Node}} & \textbf{e-MA} & \textbf{Mat-MA} & \textbf{e-$\overline{\text{MB}}$} & \textbf{Mat-$\overline{\text{MB}}$} &       & \textbf{e-MA} & \textbf{Mat-MA} & \textbf{e-$\overline{\text{MB}}$} & \textbf{Mat-$\overline{\text{MB}}$} &       & \textbf{e-MA} & \textbf{Mat-MA} & \textbf{e-$\overline{\text{MB}}$} & \textbf{Mat-$\overline{\text{MB}}$} \\
		\midrule
		\midrule
		\textbf{21} & 1     & 0.51  & 1     & 0.42  &       & 2     & 0.43  & 1     & 0.43  &       & 2     & 0.58  & 1     & 0.46 \\
		\textbf{44} & 16    & 3     & 15    & 3     &       & 17    & 3     & 13    & 2     &       & 20    & 3     & 10    & 3 \\
		\textbf{56} & 24    & 4     & 19    & 4     &       & 18    & 4     & 15    & 3     &       & 19    & 5     & 14    & 3 \\
		\textbf{72} & 36    & 11    & 29    & 8     &       & 45    & 9     & 30    & 6     &       & 64    & 11    & 36    & 7 \\
		\textbf{90} & 62    & 18    & 59    & 16    &       & 73    & 23    & 59    & 19    &       & 99    & 21    & 64    & 14 \\
		\textbf{106} & 123   & 32    & 118   & 37    &       & 151   & 36    & 111   & 28    &       & 210   & 44    & 141   & 29 \\
		\textbf{120} & 143   & 37    & 135   & 40    &       & 144   & 39    & 121   & 36    &       & 138   & 43    & 100   & 38 \\
		\textbf{163} & 445   & 183   & 392   & 102   &       & 417   & 132   & 369   & 98    &       & 636   & 167   & 392   & 104 \\
		\textbf{182} & 577   & 212   & 481   & 140   &       & 932   & 178   & 701   & 128   &       & 1326  & 172   & 1010  & 136 \\
		\textbf{203} & 956   & 206   & 879   & 198   &       & 616   & 217   & 552   & 172   &       & 515   & 210   & TL    & 194 \\
		\textbf{254} & 6103  & 414   & TL    & 462   &       & TL    & 589   & TL    & 356   &       & TL    & 568   & TL    & 514 \\
		\textbf{264} & 3903  & 479   & TL    & 652   &       & 6453  & 470   & TL    & 405   &       & TL    & 478   & TL    & 445 \\
		\textbf{275} & 7133  & 686   & 6585  & 1082  &       & TL    & 781   & TL    & 711   &       & TL    & 819   & TL    & 636 \\
		\textbf{295} & TL    & 726   & TL    & 1068  &       & TL    & 729   & TL    & 733   &       & TL    & 836   & TL    & 775 \\
		\textbf{326} & TL    & 893   & TL    & 1027  &       & TL    & 2407  & TL    & 856   &       & TL    & 1182  & TL    & 914 \\
		\textbf{355} & TL    & 1447  & TL    & 1378  &       & TL    & 1499  & TL    & 1339  &       & TL    & 1619  & TL    & 1333 \\
		\textbf{388} & TL    & 1611  & TL    & 1711  &       & TL    & 1835  & TL    & 1726  &       & TL    & 1931  & TL    & 2115 \\
		\textbf{410} & TL    & 2133  & TL    & 2251  &       & TL    & 4062  & TL    & 2167  &       & TL    & 2847  & TL    & 1901 \\
		\textbf{436} & TL    & 2439  & TL    & 2413  &       & TL    & 2531  & TL    & 2340  &       & TL    & 3608  & TL    & 2220 \\
		\textbf{449} & TL    & 2633  & TL    & 2857  &       & TL    & 2953  & TL    & 2546  &       & TL    & 2903  & TL    & 2534 \\
		\textbf{472} & TL    & 2955  & TL    & 3185  &       & TL    & 3167  & TL    & 2971  &       & TL    & 3194  & TL    & 3155 \\
		\textbf{482} & TL    & 5054  & TL    & 4044  &       & TL    & 3715  & TL    & 3216  &       & TL    & 3885  & TL    & 3162 \\
		\textbf{500} & TL    & 6238  & TL    & 4442  &       & TL    & 5653  & TL    & 2955  &       & TL    & 4278  & TL    & 2358 \\
		\bottomrule
		\bottomrule
	\end{tabular}%
}
}
	\label{tab5}%
\end{table*}%

\begin{table*}[h!]
	\centering
	\caption{Performance comparison of a combination of $\epsilon$-constraint method (e), and the proposed matheuristic (Mat) method with MA and $\overline{\text{MB}}$ for $\alpha$-level 0.7 and its corresponding $k$-value for larger number of scenarios. T[s] indicates the run time in seconds.}
	\small\addtolength{\tabcolsep}{-3pt}
	\scalebox{0.85}{
		\begin{tabular}{ccccccccccccc}
		& \multicolumn{11}{l}{$\alpha$/$k$ value: 0.7/3}                                                   &  \\
		\cmidrule{2-13}          &       & \multicolumn{2}{c}{\textbf{e-MA}} &       & \multicolumn{2}{c}{\textbf{Mat-MA}} &       & \multicolumn{2}{c}{\textbf{e-$\overline{\text{MB}}$}} &       & \multicolumn{2}{c}{\textbf{Mat-$\overline{\text{MB}}$}} \\
		\cmidrule{1-4}\cmidrule{6-7}\cmidrule{9-10}\cmidrule{12-13}    \textit{\#Node} & \textit{\#Scenario} & T[s]  & \#NDP &       & T[s]  & \#NDP &       & T[s]  & \#NDP &       & T[s]  & \#NDP \\
		\cmidrule{1-4}\cmidrule{6-7}\cmidrule{9-10}\cmidrule{12-13}    \multirow{3}[2]{*}{21} & 100   & 161   & 17    &       & 27    & 17    &       & 56    & 17    &       & 5     & 17 \\
		& 500   & 6476  & 16    &       & 348   & 16    &       & 322   & 16    &       & 52    & 16 \\
		& 1000  & TL    & 6     &       & 1668  & 17    &       & TL    & 13    &       & 282   & 17 \\
		\cmidrule{1-4}\cmidrule{6-7}\cmidrule{9-10}\cmidrule{12-13}    \multirow{3}[2]{*}{44} & 100   & 1076  & 30    &       & 178   & 30    &       & 316   & 30    &       & 34    & 30 \\
		& 500   & TL    & 5     &       & 1374  & 32    &       & TL    & 19    &       & 365   & 32 \\
		& 1000  & TL    & 3     &       & TL    & 3     &       & TL    & 7     &       & 1239  & 34 \\
		\cmidrule{1-4}\cmidrule{6-7}\cmidrule{9-10}\cmidrule{12-13}    \multirow{3}[2]{*}{56} & 100   & 2267  & 39    &       & 234   & 39    &       & 583   & 39    &       & 51    & 39 \\
		& 500   & TL    & 3     &       & 2675  & 40    &       & TL    & 22    &       & 776   & 38 \\
		& 1000  & TL    & 3     &       & TL    & 2     &       & TL    & 3     &       & 1747  & 41 \\
		\cmidrule{1-4}\cmidrule{6-7}\cmidrule{9-10}\cmidrule{12-13}    \multirow{3}[2]{*}{72} & 100   & 4186  & 50    &       & 180   & 50    &       & 826   & 50    &       & 107   & 50 \\
		& 500   & TL    & 3     &       & TL    & 2     &       & TL    & 15    &       & 983   & 51 \\
		& 1000  & TL    & 2     &       & TL    & 2     &       & TL    & 2     &       & TL    & 2 \\
		\bottomrule
		\bottomrule
	\end{tabular}%
	}
	\label{tab6}%
\end{table*}%

Performance comparison of the Mat method and the $\epsilon$-constraint method for small instances where the latter could find the optimal Pareto frontier is shown in Tables~\ref{tab8} and \ref{tab9}. To measure the quality, we report the hypervolume gap (gH\%) and the value of the $I_\epsilon$ indicator. Analyzing the results in Table~\ref{tab8} confirms that the combination of Mat and MA gives an upper bound for the Pareto frontier (as Mat is a heuristic), whereas the combination of the $\epsilon$-constraint method with $\overline{\text{MB}}$ finds a lower bound for the Pareto frontier. A lower bound is obtained because the $\epsilon$-constraint method computes exact solutions of the approximated problem, but approximates the CVaR objective from below, as the $\overline{\text{MB}}$ only works with a reduced number of subsets that are not necessarily worst for each non-dominated point. Therefore, we obtain \textit{negative} gH\% values and $I_\epsilon$ values $\leq 1$. Nonetheless, both methods provide high-quality approximations of Pareto frontiers as both gH\%, and $I_\epsilon$ represent a gap between the approximated sets and the reference set. 

\begin{table*}[!h]
	\centering
	\caption{Quality comparison using $I_\epsilon$ and hypervolume gap (gH\%) for approximated Pareto frontiers where the optimal Pareto frontier is known from solutions of solved instances using $\epsilon$-constraint method (e), with the sample size of 10 scenarios}
	\scalebox{0.85}{
		\begin{tabular}{c|cccccrcc}
			\multicolumn{1}{r}{} & \multicolumn{1}{l}{$\alpha$/$k$ value: 0.7/3} &       &       &       &       &       &       &  \\
			\cmidrule{2-9}    \multicolumn{1}{r}{} & \multicolumn{2}{l}{\textbf{e-MA (Reference set)}} &       & \multicolumn{1}{l}{\textbf{Mat-MA}} &       &       & \multicolumn{1}{l}{\textbf{e-$\overline{\text{MB}}$}} &  \\
			\cmidrule{2-3}\cmidrule{5-6}\cmidrule{8-9}    \multicolumn{1}{c}{\textit{\#Node}} & \textbf{gH\%} & \textbf{$I_\epsilon$} &       & \textbf{gH\%} & \textbf{$I_\epsilon$} &       & \textbf{gH\%} & \textbf{$I_\epsilon$} \\
			\cmidrule{1-3}\cmidrule{5-6}\cmidrule{8-9}    \textbf{21} & 0     & 1     &       & 0.00  & 1.00  &       & -0.01 & 0.99 \\
			\textbf{44} & 0     & 1     &       & 0.00  & 1.00  &       & 0.00  & 1.00 \\
			\textbf{56} & 0     & 1     &       & 0.00  & 1.00  &       & 0.00  & 0.91 \\
			\textbf{72} & 0     & 1     &       & 0.42  & 1.22  &       & -0.02 & 0.95 \\
			\textbf{90} & 0     & 1     &       & 0.05  & 1.32  &       & 0.00  & 0.87 \\
			\textbf{106} & 0     & 1     &       & 0.06  & 1.33  &       & -0.01 & 0.94 \\
			\textbf{120} & 0     & 1     &       & 0.02  & 1.35  &       & 0.00  & 1.00 \\
			\textbf{163} & 0     & 1     &       & 0.00  & 1.00  &       & 0.00  & 0.96 \\
			\textbf{182} & 0     & 1     &       & 0.00  & 1.00  &       & 0.00  & 0.94 \\
			\textbf{203} & 0     & 1     &       & 0.02  & 1.02  &       & 0.00  & 0.97 \\
			\bottomrule
			\bottomrule
		\end{tabular}%
	}
	\label{tab8}%
\end{table*}%

To measure the performance of the combination of the Mat method and $\overline{\text{MB}}$, we re-evaluate the second stage of $\overline{\text{MB}}$ by solving the MA model, where the first stage decisions are fixed to the values obtained from the $\overline{\text{MB}}$. In other words, in order to assess its true quality, the solution provided by $\overline{\text{MB}}$ is re-evaluated by using the exact second-stage objective value as obtained through MA. As can be seen in Table~\ref{tab9}, the e-$\overline{\text{MB}}$ and Mat-$\overline{\text{MB}}$ frameworks generate high-quality Pareto frontiers.

\begin{table*}[!h]
	\centering
	\caption{Quality comparison using $I_\epsilon$ and hypervolume gap (gH\%) for re-evaluated approximated Pareto frontiers using e-$\overline{\text{MB}}$ and Mat-$\overline{\text{MB}}$, where the optimal Pareto frontier is known from solutions of solved instances using $\epsilon$-constraint method (e), with the sample size of 10 scenarios}
	\scalebox{0.85}{
		\begin{tabular}{c|cccccrcc}
			\multicolumn{1}{r}{} & \multicolumn{1}{l}{$\alpha$/$k$ value: 0.7/3} &       &       &       &       &       &       &  \\
			\cmidrule{2-9}    \multicolumn{1}{r}{} & \multicolumn{2}{l}{\textbf{e-MA (Reference set)}} &       & \multicolumn{1}{l}{\textbf{e-$\overline{\text{MB}}$}} &       &       & \multicolumn{1}{l}{\textbf{Mat-$\overline{\text{MB}}$}} &  \\
			\cmidrule{2-3}\cmidrule{5-6}\cmidrule{8-9}    \multicolumn{1}{c}{\textit{\#Node}} & \textbf{gH\%} & \textbf{$I_\epsilon$} &       & \textbf{gH\%} & \textbf{$I_\epsilon$} &       & \textbf{gH\%} & \textbf{$I_\epsilon$} \\
			\cmidrule{1-3}\cmidrule{5-6}\cmidrule{8-9}    \textbf{21} & 0     & 1     &       & 0.00  & 1.00  &       & 0.00  & 1.00 \\
			\textbf{44} & 0     & 1     &       & 0.00  & 1.00  &       & 0.00  & 1.00 \\
			\textbf{56} & 0     & 1     &       & 0.00  & 1.00  &       & 0.00  & 1.00 \\
			\textbf{72} & 0     & 1     &       & 0.01  & 1.03  &       & 0.01  & 1.03 \\
			\textbf{90} & 0     & 1     &       & 0.00  & 1.04  &       & 0.00  & 1.04 \\
			\textbf{106} & 0     & 1     &       & 0.01  & 1.00  &       & 0.01  & 1.00 \\
			\textbf{120} & 0     & 1     &       & 0.00  & 1.00  &       & 0.00  & 1.00 \\
			\textbf{163} & 0     & 1     &       & 0.00  & 1.00  &       & 0.00  & 1.00 \\
			\textbf{182} & 0     & 1     &       & 0.00  & 1.04  &       & 0.00  & 1.04 \\
			\textbf{203} & 0     & 1     &       & 0.00  & 1.00  &       & 0.00  & 1.00 \\
			\bottomrule
			\bottomrule
		\end{tabular}%
	}
	\label{tab9}%
\end{table*}%

After that, we also assess the performance of Mat-$\overline{\text{MB}}$ (with re-evaluated second stage) on large-size instances where the optimal Pareto frontier is not known. For this purpose, we consider a union of the Pareto frontiers ($U$) obtained by e-MA, BB-MA, and Mat-MA as the reference set. Table~\ref{tab10} shows the hypervolume gap (gH\%) values. Negative values indicate that Mat-$\overline{\text{MB}}$ with re-evaluated second stage values obtains better approximations than the union of the other methods.

\begin{table}[!h]
	\centering
	
	\caption{Performance of Mat-$\overline{\text{MB}}$ on large-size instances, where the optimal Pareto-frontier is not known. The union of the Pareto frontiers obtained by e-MA, BB-MA, and Mat-$\overline{\text{MB}}$ is considered as a reference set. The sample size of scenarios is 10}
	\scalebox{0.85}{
		\begin{tabular}{c|ccc}
			\multicolumn{1}{r}{} & \multicolumn{1}{l}{$\alpha$/$k$ value: 0.7/3} &       &  \\
			\cmidrule{2-4}    \multicolumn{1}{r}{} & \multicolumn{1}{l}{\textbf{U (Reference set)}} &       & \multicolumn{1}{l}{\textbf{Mat-$\overline{\text{MB}}$}} \\
			\cmidrule{2-2}\cmidrule{4-4}    \multicolumn{1}{c}{\textit{\#Node}} & \textbf{gH\%} &       & \textbf{gH\%} \\
			\cmidrule{1-2}\cmidrule{4-4}    \textbf{254} & 0.00  &       & 0.00 \\
			\textbf{275} & 0.00  &       & 0.00 \\
			\textbf{295} & 0.00  &       & -0.02 \\
			\textbf{326} & 0.00  &       & -0.06 \\
			\textbf{355} & 0.00  &       & 0.01 \\
			\textbf{388} & 0.00  &       & -0.10 \\
			\textbf{410} & 0.00  &       &  -0.05\\
			\textbf{436} & 0.00  &       &  -0.02\\
			\textbf{449} & 0.00  &       & -0.03 \\
			\textbf{472} & 0.00  &       & -0.08 \\
			\textbf{482} & 0.00  &       & -0.02 \\
			\textbf{500} & 0.00  &       &  0.00\\
			\bottomrule
			\bottomrule
		\end{tabular}%
		
	}
	\label{tab10}%
\end{table}%

\section{\uppercase{Conclusions}}
\label{sec:conclusion}

In this paper, we investigate the subset-based polyhedral representation of CVaR to capture risk in a two-stage bi-objective location-allocation model introduced in \cite{nazemi2020bi}, the purpose of which is the optimal design of a relief network. The model aims at finding the trade-off between a deterministic objective (minimization of location cost) and an uncertain one (minimization of uncovered demand) where uncertainty is assumed in demand values. We introduce an approximate cutting-plane method to deal with the computationally challenging subset-based formulation of the two-stage stochastic optimization problem containing the CVaR.  

We evaluate the performance of the proposed method by conducting a computational study on the instances used by \cite{nazemi2020bi}. Our experiments show that the proposed approximate cutting-plane approach pays off for a large number of scenarios and they illustrate how incorporating an approximation of the subset-based formulation affects the quality of the produced solutions. In future research, we would like to investigate whether the performance of the cutting-plane method can be further improved by employing additional enhancements.

\subsection*{\uppercase{Acknowledgements}}
The authors want to thank \cite{tricoire2012bi} for providing us with the real-world data set. This research was funded in whole, or in part, by the Austrian Science Fund (FWF) [P 31366]. For the purpose of open access, the author has applied a CC BY public copyright licence to any Author Accepted Manuscript version arising from this submission.

 \bibliographystyle{apalike}
 {\small

 \bibliography{ICORESpaperRef}
 }

\end{document}